\newif\ifarXiv    
\newif\ifWP       
\arXivtrue

\newif\ifTR     
\ifarXiv\TRtrue\fi
\ifWP\TRtrue\fi

\newif\ifFULL

\ifarXiv
  \documentclass[10pt]{article}
\fi

\ifWP
  \documentclass[10pt]{article}
  \usepackage{amsmath,amsthm,amsfonts,amssymb,latexsym,graphicx,epigraph,accents}
  \usepackage{etoolbox}
  \input{/Doc/Computing/Latex/kp.txt}  
\fi

\ifTR
  \usepackage{amsmath,amsfonts,amssymb,amsthm,graphicx}
  \usepackage{enumerate}
  \usepackage{etoolbox}
  \usepackage{color}
  \usepackage[colorlinks=true,citecolor=blue,pdfpagemode=UseNone,pdfstartview=FitH]{hyperref}
\fi

\newtoggle{arXiv}
\newtoggle{WP}
\newtoggle{FULL}
\newtoggle{TR}
\newtoggle{false} 
\newtoggle{true}\toggletrue{true} 

\ifarXiv\toggletrue{arXiv}\fi
\ifWP\toggletrue{WP}\fi
\ifFULL\toggletrue{FULL}\fi
\ifTR\toggletrue{TR}\fi

\iftoggle{FULL}{
  \newcommand{\bluebegin}{\begingroup\color{blue}}
  \newcommand{\blueend}{\endgroup}

}{}

\allowdisplaybreaks[4]

\newcommand{\R}{\mathbb{R}}    

\iftoggle{TR}{%
  \theoremstyle{plain}
  \newtheorem{theorem}{Theorem}[section]

  \theoremstyle{definition}
  
  \newtheorem{example}[theorem]{Example}
  \theoremstyle{remark}
  
}{}

\newcommand{\abstr}%
  {A very simple example demonstrates
  that Fisher's application of the conditionality principle to regression
  (``fixed-$x$ regression''),
  endorsed by David Sprott and many other followers,
  makes prediction impossible in the context of statistical learning theory.
  On the other hand, relaxing the requirement of conditionality
  makes it possible
  via, e.g., conformal prediction.}

\iftoggle{WP}{%
  \twodatestrue
  
}{}

\begin{document}
\iftoggle{arXiv}{%
  \title{Conditionality principle under unconstrained randomness}
  \author{Vladimir Vovk}
  \maketitle
  \begin{abstract}
    \smallskip
    \abstr

    The version of this note at \url{http://alrw.net} (Working Paper 40)
    is updated most often.
  \end{abstract}
}{}

\iftoggle{WP}{%
  \title{Conditionality principle under unconstrained randomness}
  \author{Vladimir Vovk}
  \maketitle
  \begin{abstract}
    \smallskip
    \abstr
  \end{abstract}
}{}

\section{Introduction}

The main goal of this note is to draw the reader's attention
to the fact that the conditionality principle
is not compatible with statistical learning theory,
in which the task is to predict the label $y$ of an object $x$.
Two characteristic features of statistical learning theory
are that the labelled objects $(x,y)$ are only assumed
to be independent and identically distributed
(the unrestricted \emph{assumption of randomness})
and that the objects $x$ are complex (such as videos),
so that we are unlikely to ever see identical objects.
These two features make prediction impossible
if we want to condition on the observed sequence of $x$s as recommended by Fisher.
This is not a new observation\iftoggle{FULL}{\bluebegin\
  (it has been known for a long time)\blueend}{},
but it might not be as widely known
as it deserves.

\iftoggle{FULL}{\bluebegin
  I start in Sect.~\ref{sec:IID} by discussing
  the standard framework of machine learning and the assumption of randomness.
  A simple conditionality principle is stated and discussed in Sect.~\ref{sec:principle-1}.
  In the same section I give an example highlighting the impossibility of prediction.
  The example is very simple
  and does not rely on any previous knowledge of conformal prediction.
  Section~\ref{sec:results} summarizes,
  and gives references to,
  some of the numerous technical results in this area.
  Finally, in Sect.~\ref{sec:conclusion} I put forward my suggestions
  for resolving the contradiction.%
\blueend}{}%

\section{Assumption of randomness and conformal prediction}
\label{sec:IID}

In statistical learning theory
(see, e.g., \cite{Vapnik/Chervonenkis:1974,Vapnik:1998})
we consider \emph{observations} $(x,y)$ each consisting of two components:
an \emph{object} $x\in\mathbf{X}$ and its \emph{label} $y\in\mathbf{Y}$.
In general, the \emph{object space} $\mathbf{X}$ and \emph{label space} $\mathbf{Y}$
are arbitrary measurable spaces,
but to discuss the relevance to Fisher's ideas
it will often be convenient to concentrate on the case of \emph{regression} $\mathbf{Y}=\R$.

The simplest setting
is where we are given
a training sequence
\[
  (x_1,y_1),\dots,(x_n,y_n)%
\]
and the problem is to predict the label $y_{n+1}$ of a test object $x_{n+1}$.
The (unrestricted) \emph{assumption of randomness}
is that the observations $(x_1,y_1),\dots,(x_{n+1},y_{n+1})$
are generated independently from an unknown probability measure on $\mathbf{X}\times\mathbf{Y}$.
This assumption is standard in machine learning
and popular in nonparametric statistics.

One way to make predictions
with validity guarantees under unconstrained randomness
is \emph{conformal prediction} \cite{Angelopoulos/etal:2024-local,Vovk/etal:2022book}:
given a target probability of error $\epsilon>0$
conformal prediction produces a prediction set $\Gamma\subseteq\mathbf{Y}$
such that $y_{n+1}\in\Gamma$ with probability at least $1-\epsilon$.
The basic idea of conformal prediction is familiar
(see, e.g., \cite[Sect.~7.5]{Cox/Hinkley:1974}):
we fix a statistical test of the null hypothesis of randomness,
go over all possible labels $y$ for the test object $x_{n+1}$,
and include in $\Gamma$ all labels $y$
for which the augmented training set
$(x_1,y_1),\dots,(x_n,y_n),(x_{n+1},y)$
does not lead to the rejection of the null hypothesis.
In many interesting cases this idea has
computationally efficient implementations
(see, e.g., \cite[Sects.~2.3 and~2.4]{Vovk/etal:2022book},
\cite[Sect.~9.2]{Angelopoulos/etal:2024-local},
and \cite{Lei:2019}).

\iftoggle{FULL}{\bluebegin
  Other ways to make predictions
  with validity guarantees under unconstrained randomness
  include numerous results in statistical learning theory and PAC theory,
  starting from Vapnik and Chervonenkis.

  It is interesting that Fisher was downplaying nonparametric statistics
  even when introducing pioneering nonparametric methods,
  as in Sects.~21 and~21.1 of his 1935 book \cite{Fisher:1935book}.%
\blueend}{}%

\section{A simple conditionality principle and example}
\iftoggle{FULL}{\label{sec:principle-1}}{}%

In this section we will only need
a special case of the conditionality principle,
which I will call the ``fixed-$x$ principle''
partly following Aldrich \cite{Aldrich:2005}.
The fixed-$x$ principle says that,
when performing any kind of statistical analysis,
we should consider the observed sequence of objects
$x_1,\dots,x_{n+1}$
as fixed,
even if they were in fact generated from some probability distribution
(known or unknown).
This is applicable to regression problems (\emph{fixed-$x$ regression}),
or any other prediction problems of the kind
described in the previous section.

Fisher was a life-long promoter
of both the general conditionality principle
(which was introduced formally only in 1962 by Birnbaum \cite{Birnbaum:1962})
and the fixed-$x$ principle as its special case.
\iftoggle{FULL}{\bluebegin
  Even though the fixed-$x$ principle was widely used by Fisher
  and his disciples (such as Sprott),
  it is very difficult to pinpoint where he stated it clearly.
\blueend}{}%
When H. Fairfield Smith asked Fisher
about the origin of the fixed-$x$ principle
(treating ``the independent variable as fixed
even although it might have been observed
as a random sample of some variate population'')
in his letter of 6 August 1954 \cite[pp.~213--214]{Bennett:1990},
Fisher responded with a reference to his 1922 paper
\cite[p.~599]{Fisher:1922-local}.

\iftoggle{FULL}{\bluebegin
  This is the quote from Smith (which is too long to give in full):
  \begin{quote}
    I have long been under the impression
    that when I was at the Galton Laboratory
    (and not yet able to appreciate the points involved),
    I saw a short paper of two or three pages by you
    claiming that for purposes of a regression analysis
    it would always be valid to treat the independent variable as fixed
    even although it might have been observed
    as a random sample of some variate population.
    The point is frequently cropping up,
    although not usually of much importance,
    and I am unable to quote your reference if it exists
    or to track it down in any bibliography.
    I have always meant to ask you when we meet but invariably forget.
    I have asked many other statisticians and none know the reference,
    although you are frequently quoted as having given the conclusion
    as a verbal opinion.\dots
  \end{quote}
  (See \cite[pp.~213--214]{Bennett:1990}.)
  Fisher referred Smith to his 1922 paper \cite[p.~599]{Fisher:1922-local}.
  (But what Fisher said is slightly obscure.)%
\blueend}{}%

See Aldrich \cite{Aldrich:2005} about the development of Fisher's fixed-$x$ regression.
Aldrich quotes Fisher \cite[Sect.~2, p.~71]{Fisher:1955}
(in the context of regression with $y$ having a Gaussian distribution given $x$):
``The qualitative data may also tell us how $x$ is distributed,
with or without specific parameters;
this information is irrelevant.''

David Sprott, a prominent follower of Fisher's,
also promoted the conditionality principle in his work.
In his 1989 interview with Mary Thompson \cite{Thompson:1989},
he remembers a case when, as a student,
he was tempted to take into account the variation in the $x$s
in a practical regression problem.
His statistics professor said,
``No, you wouldn't do that, you'd condition on the $x$s''.
Sprott couldn't find out why conditioning on the $x$s was the right thing to do
until he went to London a few years later to work with Fisher,
but then he was fully convinced by Fisher.

\iftoggle{FULL}{\bluebegin
  This is the full quote from Sprott in his interview:
  \begin{quote}
    ``\dots I do remember a case when in the fourth year somebody brought me
    a practical problem in regression in which the $x$'s were uncontrolled,
    and it had occurred to me, following rigorously the theory I had been taught,
    that you should allow for the variation in the $x$'s
    by taking the joint distribution of the $x$'s and $y$'s;
    and the statistics professor said
    ``No, you wouldn't do that, you'd condition on the $x$'s,''
    and I said, ``But why? The course you've given us wouldn't say you should do that at all.''
    And he said ``No, but that's what you'd do.''
    Well that's what Fisher said too, but he was saying why you would do it.''
  \end{quote}
\blueend}{}%

The intuition behind the fixed-$x$ principle
is that only the observed objects are relevant for predicting the label of $x_{n+1}$.
In Sect.~\ref{sec:principle-2} we will discuss this in a wider context.
And in this section, we discuss
the paralysing effect of the fixed-$x$ principle under unrestricted randomness
using the following example.

\begin{example}\label{ex:example}
  Consider the problem of regression with $\mathbf{X}=\mathbf{Y}=\R$,
  and suppose:
  \begin{itemize}
  \item
    the training sequence is such that $y_i=x_i$ for all $i=1,\dots,n$, for a large~$n$;
  \item
    $x_1,\dots,x_{n+1}$ are all different.
  \end{itemize}
  What can we say about $y_{n+1}$ knowing $x_{n+1}$?
\end{example}

Example~\ref{ex:example} may describe a situation
where the observations are independent and coming from the same continuous distribution.
Under the assumption of randomness,
we can confidently claim that $y_{n+1}=x_{n+1}$;
otherwise, the last observation looks strange
and leads to a p-value of $1/(n+1)$ for a fixed statistical test.
As the test statistic $T$ for such a test we can take, e.g.,
\[
  T
  :=
  \begin{cases}
    1 & \text{if $\lvert y_{n+1}-x_{n+1}\rvert>\lvert y_i-x_i\rvert$, $i=1,\dots,n$,}\\
    0 & \text{otherwise}.
  \end{cases}
\]
For any $\epsilon\ge1/(n+1)$,
we have $[x_{n+1},x_{n+1}]$ as the prediction interval for $y_{n+1}$
at confidence level $1-\epsilon$.
Intuitively, this follows from the expectation that the future will be similar to the past
(Laplace's rule of succession).
Conformal prediction extends this idea greatly.

However, we can say nothing whatsoever about $y_{n+1}$
if we condition on the observed $x_1,\dots,x_{n+1}$.
The problem with conditioning on $x$ is that it destroys the assumption of randomness.
The assumption of randomness becomes the following \emph{assumption of conditional randomness}:
the data-generating distribution $P$ is determined by an arbitrary sequence
$(x_1,\dots,x_{n+1})\in\mathbf{X}^{n+1}$
and an arbitrary family of probability measures
$\{Q_x\mid x\in\mathbf{X}\}$ on $\mathbf{Y}$
(measurable in the sense of being a Markov kernel
with $\mathbf{X}$ as source and $\mathbf{Y}$ as target);
we then have
$P=(\delta_{x_1}\times Q_{x_1})\times\dots\times(\delta_{x_{n+1}}\times Q_{x_{n+1}})$,
where $\delta_x$ is the distribution on $\mathbf{X}$
that is concentrated at $x$.
When $x_1,\dots,x_{n+1}$ are all different,
the true conditional distribution of $y_{n+1}$
can be any probability measure on $\mathbf{Y}$
as there are no restrictions on $Q_{x_{n+1}}$.


\section{Different kinds of conditionality}
\iftoggle{FULL}{\label{sec:results}}{}%

In this section I will briefly describe a couple of results that shed light
on the phenomenon illustrated
by Example~\ref{ex:example}.
But let me first embed the fixed-$x$ principle in a more general picture.

In conformal prediction,
a common goal is to achieve \emph{conditional validity},
i.e., to make the conditional probability of an error,
$y_{n+1}\notin\Gamma$ in the notation of Sect.~\ref{sec:IID},
bounded above by a given $\epsilon$.
But what should we condition on?
Three kinds of conditional validity have been widely studied
in conformal prediction
(see, e.g., \cite[Sect.~4.7.1]{Vovk/etal:2022book} or,
in greater detail, \cite[Chap.~4]{Angelopoulos/etal:2024-local}):
\begin{itemize}
\item
  in \emph{object-conditional validity},
  we condition on $x_{n+1}$
  (or on a function thereof,
  in which case we will talk about \emph{partial object-conditional validity});
\item
  in \emph{label-conditional validity},
  we condition on $y_{n+1}$;
\item
  in \emph{training-conditional validity},
  we condition on $x_1,y_1,\dots,x_n,y_n$.
\end{itemize}
It might be possible to achieve more than one kind of conditional validity at the same time.
One such case is the Bayesian setting where, e.g.,
in addition to the assumption of randomness
we postulate a probability measure over the probability measures $R$
generating one observation;
under this extra assumption, we can construct efficient prediction sets
satisfying training-conditional and object-conditional validity at the same time.
But of course, it does not make sense to aim at achieving all three kinds of conditional validity
at the same time.

In this note we are interested in a fourth kind of conditional validity,
in which we condition on $x_1,\dots,x_{n+1}$.
This is a stronger requirement than object-conditional validity,
and already the latter is difficult to achieve under unrestricted randomness.

Lei and Wasserman
(\cite[Lemma~1]{Lei/Wasserman:2014}, \cite[Theorem 4.11]{Vovk/etal:2022book})
prove a strong negative result assuming that $\mathbf{X}$ is a separable metric space.
In the case of regression ($\mathbf{Y}=\R$),
they show that a set predictor satisfying
the property of object-conditional validity under unrestricted randomness
will output, with probability at least $1-\epsilon$,
an unbounded prediction set for $y_{n+1}$
provided the test object $x_{n+1}$ is not an atom of the data-generating distribution.
Remember that $\epsilon$ is the target probability of error,
so that the lower bound $1-\epsilon$ means complete lack of efficiency.

Barber et al.\ \cite{Barber/etal:2021limits} show
that an efficient set predictor can satisfy partial object-conditional validity
that involves conditioning on $x_{n+1}\in\mathcal{X}$
for all sets $\mathcal{X}$ of probability at least $\delta$
for a given $\delta>0$,
but essentially the same properties of efficiency and validity
are achieved automatically by unconditionally valid predictors,
such as conformal predictors
(under some regularity conditions).

Despite these negative results,
designing conformal predictors that are conditional in a weaker sense
is an active area of research,
starting from a basic idea of Mondrian conformal prediction
(\cite[Sect.~4.6]{Vovk/etal:2022book},
\cite[Chap.~4]{Angelopoulos/etal:2024-local}).
Asymptotically, conditional conformal prediction
is possible in a very strong sense \cite[Theorem~1]{Lei/Wasserman:2014}.
Finite-sample results are mathematically less satisfying
but may hold great promise in practice;
see, e.g., \cite{Romano/etal:2020,Bellotti:2021,Guan:2023,Hore/Barber:arXiv1904}.

\section{General conditionality principle and its difficulties}
\label{sec:principle-2}

\iftoggle{FULL}{\bluebegin
  According to Wikipedia,
  the conditionality principle is a Fisherian principle of statistical inference.%
\blueend}{}%

Cox and Hinkley \cite[Sect.~2.3(iii)]{Cox/Hinkley:1974}
point out two versions of the conditionality principle,
basic and extended.
In the basic version,
we are given an \emph{ancillary statistic} $C$,
i.e., a random variable with a known distribution,
and the conditionality principle says that our analysis
should be conditional on the observed value of $C$.
This prescription is very compelling in some cases,
such as Cox's famous example
of choosing one of two measuring instruments at random
and then observing its reading
(knowing which instrument has been chosen);
see \cite[Example~2.33]{Cox/Hinkley:1974}.

The basic version does not imply the fixed-$x$ principle,
since the distribution of the $x$s
does depend on the unknown data-generating distribution.
In the extended conditionality principle,
the unknown parameter is split into two parts,
and only one of those parts is of direct interest to us.
In the context of the assumption of randomness,
the parameter is the probability measure $R$ on $\mathbf{X}\times\mathbf{Y}$
generating one observation.
We can typically (see, e.g., \cite[Sect.~A.4]{Vovk/etal:2022book}) split $R$
into the marginal distribution $R_{\mathbf{X}}$ on $\mathbf{X}$
and the conditional distribution $Q_x$ of $y$ given $x$, for each $x\in\mathbf{X}$;
$Q$ is a Markov kernel.
Only $Q:=\{Q_x\mid x\in\mathbf{X}\}$ is of interest to us in our prediction problem.
Then $C:=(x_1,\dots,x_{n+1})$ is ancillary for the Markov kernel $Q$
in the following extended sense:
\begin{itemize}
\item
  The distribution of $C$ does not depend on $Q$
  (and only depends on $R_{\mathbf{X}}$).
\item
  The conditional distribution of the remaining part of the data
  $y_1,\dots,y_{n+1}$ given the value of $C$
  depends only on $Q$
  (and does not depend on $R_{\mathbf{X}}$).
\end{itemize}
The extended conditionality principle,
requiring analysis to be conditional on $C$,
becomes the fixed-$x$ principle
when adapted to the assumption of randomness.

\iftoggle{FULL}{\bluebegin
  Sprott (2000):
  ``Conditioning enforces the criterion of relevance''.%
\blueend}{}%

Lehmann and Scholz \cite[Sect.~1]{Lehmann/Scholz:1992}
point out that conditional inference can be less efficient for small samples,
although the difference tends to disappear as the sample size increases.
However, Example~\ref{ex:example} is more serious,
since it demonstrates a complete failure of the conditionality principle.
Another instance of a comparable complete failure of this principle
is where the experimental design involves deliberate randomization,
as in a random assignment of subjects to treatments in randomized clinical trials
\cite[end of Sect.~3]{Lehmann/Scholz:1992}.
The conditionality principle then forces us to disregard randomization
compromising one of the most standard and powerful statistical tools in medicine.

\section{Conclusion}
\iftoggle{FULL}{\label{sec:conclusion}}{}%

This note observes that the extended conditionality principle
prevents successful prediction in statistical learning theory,
which is the basic setting of machine learning.
(Other varieties of machine learning usually make prediction even more difficult;
e.g., they may allow different distributions for the training and test observations.)
Its main points (which are far from being original) are:
\begin{itemize}
\item
  we should drop the conditionality principle as a requirement but keep it as an ideal goal
  (which is often partially achievable but more rarely,
  and then under strong assumptions,
  fully achievable);
\item
  doing so permits use of methods with unconditional performance guarantees,
  such as conformal prediction
  (and these methods automatically satisfy
  some weak properties of conditional validity \cite{Barber/etal:2021limits});
\item
  when trying to achieve conditionality as ideal goal,
  we might need to settle for partial, approximate, or asymptotic conditionality;
\item
  it is important to keep in mind other kinds of conditionality,
  such as training-conditional, object-conditional, and label-conditional validity.
\end{itemize}
For several results in these directions,
see, e.g., \cite[Chap.~4]{Angelopoulos/etal:2024-local}
and \cite[Sects.~1.4.4, 4.6, 4.7]{Vovk/etal:2022book}.

\subsection*{Acknowledgments}

Thanks to Changbao Wu,
Chair of the Department of Statistics and Actuarial Science
(University of Waterloo, Canada),
for inviting me to deliver the David Sprott Distinguished Lecture in June 2023.
I am grateful to my listeners, especially Ruodu Wang,
for their perceptive comments.
Many thanks to two anonymous referees
for their suggestions about presentation
in the journal version of this note.
My research has been partially supported by Mitie.

\iftoggle{TR}{%
  \bibliographystyle{plain}
}{}
\bibliography{local,%
  /doc/work/r/bib/general/general,%
  /doc/work/r/bib/math/math,%
  /doc/work/r/bib/prob/prob,%
  /doc/work/r/bib/stat/stat,%
  /doc/work/r/bib/vovk/vovk}

\begin{thebibliography}{10}

\bibitem{Aldrich:2005}
John Aldrich.
\newblock Fisher and regression.
\newblock {\em Statistical Science}, 20:401--417, 2005.
\newblock Original title: ``The origins of fixed $X$ regression''.

\bibitem{Angelopoulos/etal:2024-local}
Anastasios~N. Angelopoulos, Rina~Foygel Barber, and Stephen Bates.
\newblock Theoretical foundations of conformal prediction.
\newblock Technical Report \href{https://arxiv.org/abs/2411.11824}{2411.11824
  [math.ST]}, \href{https://arxiv.org/}{arXiv.org} e-Print archive, November
  2024.
\newblock Pre-publication version of a monograph to be published by Cambridge
  University Press.

\bibitem{Barber/etal:2021limits}
Rina~Foygel Barber, Emmanuel~J. Cand\`es, Aaditya Ramdas, and Ryan~J.
  Tibshirani.
\newblock The limits of distribution-free conditional predictive inference.
\newblock {\em Information and Inference: A Journal of the IMA}, 10:455--482,
  2021.

\bibitem{Bellotti:2021}
Anthony Bellotti.
\newblock Approximation to object conditional validity with inductive conformal
  predictors.
\newblock {\em Proceedings of Machine Learning Research}, 152:4--23, 2021.
\newblock {COPA} 2021.

\bibitem{Bennett:1990}
J.~H. Bennett, editor.
\newblock {\em Statistical Inference and Analysis: Selected Correspondence of
  R.~A.~Fisher}.
\newblock Clarendon Press, Oxford, 1990.

\bibitem{Birnbaum:1962}
Allan Birnbaum.
\newblock On the foundations of statistical inference (with discussion).
\newblock {\em Journal of the American Statistical Association}, 57:269--326,
  1962.

\bibitem{Cox/Hinkley:1974}
David~R. Cox and David~V. Hinkley.
\newblock {\em Theoretical Statistics}.
\newblock Chapman and Hall, London, 1974.

\bibitem{Fisher:1922-local}
Ronald~A. Fisher.
\newblock The goodness of fit of regression formulae, and the distribution of
  regression coefficients.
\newblock {\em Journal of the Royal Statistical Society}, 85:597--612, 1922.

\bibitem{Fisher:1955}
Ronald~A. Fisher.
\newblock Statistical methods and scientific induction.
\newblock {\em Journal of the Royal Statistical Society B}, 17:69--78, 1955.

\bibitem{Guan:2023}
Leying Guan.
\newblock Localized conformal prediction: a generalized inference framework for
  conformal prediction.
\newblock {\em Biometrika}, 110:33--50, 2023.

\bibitem{Hore/Barber:arXiv1904}
Rohan Hore and Rina~Foygel Barber.
\newblock Conformal prediction with local weights: randomization enables robust
  guarantees.
\newblock Technical Report \href{https://arxiv.org/abs/2310.07850}{2310.07850
  [stat.ME]}, \href{https://arxiv.org/}{arXiv.org} e-Print archive, October
  2023.
\newblock Forthcoming in the \emph{Journal of the Royal Statistical Society B}.

\bibitem{Lehmann/Scholz:1992}
Erich~L. Lehmann and F.~W. Scholz.
\newblock Ancillarity.
\newblock In Malay Ghosh and Pramod~K. Pathak, editors, {\em Current Issues in
  Statistical Inference: Essays in Honor of D. Basu}, volume~17 of {\em Lecture
  Notes---Monograph Series}, pages 32--51. Institute of Mathematical
  Statistics, Hayward, CA, 1992.

\bibitem{Lei:2019}
Jing Lei.
\newblock Fast exact conformalization of lasso using piecewise linear homotopy.
\newblock {\em Biometrika}, 106:749--764, 2019.

\bibitem{Lei/Wasserman:2014}
Jing Lei and Larry Wasserman.
\newblock Distribution-free prediction bands for non-parametric regression.
\newblock {\em Journal of the Royal Statistical Society B}, 76:71--96, 2014.

\bibitem{Romano/etal:2020}
Yaniv Romano, Matteo Sesia, and Emmanuel~J. Cand\`es.
\newblock Classification with valid and adaptive coverage.
\newblock In {\em Advances in Neural Information Processing Systems 33 (NeurIPS
  2020)}, 2020.

\bibitem{Thompson:1989}
Mary Thompson.
\newblock A conversation with {David A. Sprott}.
\newblock {\em Liaison}, 3(2), Feb\-ru\-ary 1989.
\newblock Available (in March 2025) on
  \href{https://ssc.ca/en/profile/conversation-david-sprott}{the web}.

\bibitem{Vapnik:1998}
Vladimir~N. Vapnik.
\newblock {\em Statistical Learning Theory}.
\newblock Wiley, New York, 1998.

\bibitem{Vapnik/Chervonenkis:1974}
Vladimir~N. Vapnik and Alexey~Y. Chervonenkis.
\newblock \emph{Theory of Pattern Recognition} (in Russian).
\newblock Nauka, Moscow, 1974.
\newblock {German} translation: \emph{Theorie der Zeichenerkennung}, Akademie,
  Berlin, 1979.

\bibitem{Vovk/etal:2022book}
Vladimir Vovk, Alex Gammerman, and Glenn Shafer.
\newblock {\em Algorithmic Learning in a Random World}.
\newblock Springer, Cham, second edition, 2022.

\end{thebibliography}
\end{document}